\newif
\ifconfver
 \confvertrue        

\ifconfver
    \documentclass[10pt,journal,final,twoside]{IEEEtran}
\else
    \documentclass[11pt,draftcls,onecolumn]{IEEEtran}
\fi
\usepackage{xspace,empheq,fancybox,amsmath,bbm,amssymb,epsfig,syntonly,times,amsthm, cases} 
\usepackage{psfrag,color,bm,array,cite, soul}
\usepackage{url,cite,footnote,xspace,syntonly,algorithm,algpseudocode}
\usepackage{verbatim,multirow,slashbox,enumitem}
\usepackage{algorithm}
\usepackage[T1]{fontenc}
\usepackage[dvipsnames]{xcolor}
\usepackage{graphicx}

\usepackage{graphicx}
\usepackage{epstopdf}
\usepackage{mysymbol}
\usepackage{soul}

\def\bbepsilon{{\mbox{\boldmath $\epsilon$}}}
 
\def\bbgamma{{\mbox{\boldmath $\gamma$}}}

\newcommand{\utwi}[1]{\mbox{\boldmath $#1$}}

\newcommand{\diag}{\mathsf{diag}}

\renewcommand{\hat}{\widehat}
\renewcommand{\tilde}{\widetilde}

\newcommand{\cN}{{\cal N}}

\newcommand{\cE}{{\cal E}}

\newcommand{\cM}{{\cal M}}

\newcommand{\bn}{{\bf n}}

\newcommand{\bs}{{\bf s}}
\newcommand{\bx}{{\bf x}}
\newcommand{\bu}{{\bf u}}
\newcommand{\bv}{{\bf v}}
\newcommand{\bw}{{\bf w}}
\newcommand{\bi}{{\bf i}}
\newcommand{\bz}{{\bf z}}
\newcommand{\by}{{\bf y}}

\newcommand{\bG}{{\bf G}}
\newcommand{\bJ}{{\bf J}}

\newcommand{\bH}{{\bf H}}

\newcommand{\bM}{{\bf M}}
\newcommand{\bMY}{{\bf M^Y}}
\newcommand{\bMDelta}{{\bf M^\Delta}}

\newcommand{\bQ}{{\bf Q}}

\newcommand{\bY}{{\bf Y}}

\newcommand{\bepsilon}{{\utwi{\epsilon}}}

\newcommand{\comps}{\mathbb{C}}

\newcommand{\sfT}{\textsf{T}}

\DeclarePairedDelimiterX{\norm}[1]{\lVert}{\rVert}{#1}

\usepackage{epstopdf}

\theoremstyle{definition}
\newtheorem{remark}{Remark}

\allowdisplaybreaks[4]

\begin{document}
\title{Dynamic Distribution State Estimation \\ Using Synchrophasor Data}
\author{
\IEEEauthorblockN{Jianhan Song},
	\IEEEauthorblockN{Emiliano Dall'Anese},
	\IEEEauthorblockN{Andrea Simonetto}, and	
	\IEEEauthorblockN{Hao Zhu}
    \thanks{J. Song and H. Zhu are with the Department of Electrical and Computer Engineering at The University of Texas at Austin. E. Dall'Anese is with the Department of Electrical, Computer, and Energy Engineering at the University of Colorado Boulder. A. Simonetto is with IBM Research Ireland. The work of E. Dall'Anese was supported in part by the National Renewable Energy Laboratory through the grant APUP UGA-0-41026-109. H. Zhu was partially supported by NSF ECCS-1802319.}}

\markboth{}{Song \MakeLowercase{\textit{et al.}}: Dynamic Distribution State Estimation Using Synchrophasor Data}

\renewcommand{\thepage}{}
 \maketitle
\pagenumbering{arabic}

%
\begin{abstract}

The increasing deployment of distribution-level phasor measurement units (PMUs) calls for dynamic distribution state estimation (DDSE) approaches that tap into high-rate measurements to maintain a comprehensive view of the distribution-system state in real time.
Accordingly, this paper explores the development of a fast algorithmic framework  by casting the DDSE task within the \textit{time-varying} optimization realm. The time-varying formulation involves a time-varying robustified least-squares approach, and it naturally models optimal trajectories for the estimated states under streaming of measurements. The formulation is based on a linear surrogate of the AC power-flow equations,  and it includes an element of robustness with respect to measurement outliers. The paper then leverages a first-order prediction-correction method to achieve simple online updates that can provably track the state variables from heterogeneous measurements. This online algorithm is computationally efficient as it relies on the Hessian of the cost function without computing matrix-inverse. Convergence and bounds on the estimation errors of proposed algorithm can be analytically established. 

\end{abstract}

\section{Introduction}
\label{sec:intro}

Recently, power distribution networks have witnessed an increasing connection of renewable energy sources, electric vehicles, energy storage systems, among other distributed energy resources. These transformations have propelled the development and deployment of advanced sensing, communications, and control technologies. In particular, the synchrophasor technology based on 
distribution-level phasor measurement units (PMUs) \cite{Meier} has equipped distribution systems operators with synchronized, low latency and high-resolution measurements that can be collected on a fast time scale. It is timely to design efficient and effective distribution situational awareness modules that can leverage the availability of PMU data therein. 

Distribution dynamic state estimation (DDSE) is a  fundamental tool that can enable general distribution operations and control tasks. This is similar to the operational paradigm for wide-area transmission grids \cite{PMUGomez}. Nonetheless, distribution networks have unique characteristics compared to transmission systems, such as unbalanced loads and higher resistance-to-reactance ratios \cite{primadianto2017review}. Specifically, the unbalanced nature of distribution systems necessitates the multi-phase modeling at increased dimension and coupling. Moreover, nodes with zero injections commonly exist while historic load data are often included to improve redundancy. The high variation of accuracy from both types of data can lead to numerical conditioning issues, along with the higher resistance-to-reactance ratios. Thus, the dynamic SE approaches traditionally developed for transmission grids may not be directly applicable for distribution systems with PMU data.

Traditionally, distribution SE research has been limited to the static setting, constrained by meter availability and data rates. Several efforts have been focused on addressing the numerical conditioning issue using e.g., branch current formulation \cite{baran2}, conversion to current measurements \cite{lu,baran}, zero-injection information \cite{lin2,EqConKorres}. More recent work has considered the incorporation of synchrophasor data for static SE, but using the purely linear SE assuming sufficient observability from PMU data only \cite{HaughtonHeydt,LPFMeier,DSSEmuscas}. This is not yet a reality for distribution networks, and thus legacy meters should still be included for distribution SE. 
Motivated by the availability of fast PMU measurements, dynamic SE methods has been developed relying on the recursive Kalman filter updates for distribution systems \cite{sarri2012state,carquex2018state}. Note that dynamic SE in transmission systems also uses variations of Kalman filtering recursions as the ``workhorse'' algorithm; see e.g., \cite{zhao2017robust,valverde2011unscented}. Please see a recent review on distribution SE in  \cite{primadianto2017review}.

This paper aims to develop an efficient and effective DDSE solution technique that can address the unique characteristics of distribution networks and fast sampling rates of PMU data. We leverage the linearized multi-phase AC power flow model recently developed in~\cite{bernstein2017load}. The model can account for unbalanced operation, as well as for wye and delta connections. Based on this model, we formulate a new DDSE problem by advocating a \emph{time-varying} optimization formalism. The time-varying formulation involves a time-varying robustified least-squares approach, and it naturally models optimal trajectories for the estimated states under streaming of measurements.
Based on the proposed time-varying optimization model for the DDSE task, this paper proposes an algorithmic framework to \emph{track} the state of the distribution system by leveraging running prediction-correction methodologies~\cite{Simonetto_pc15,Simonetto_pc17}. Prediction-correction methods involves two phases, sequentially implemented at each time step: i) a \emph{prediction} phase where, based on the measurements collected up the current time instant, the algorithm attempts to predict the optimal  solution of the next time period  by exploring into intrinsic temporal correlations of the cost function; and, ii) once a new datum/measurement becomes available, the \emph{correction} phase refines the predicted solution. To facilitate the development of computationally affordable algorithms, the paper considers first-order prediction-correction (FOPC) methods that rely on the Hessian of the cost function, instead of requiring the computation of its inverse~\cite{Simonetto_pc17}. FOPC methods are very attractive for DDSE problems where measurements are collected at high frequency by distribution-level PMUs or even distributed energy resources to enhance real-time situational awareness. 
In particular, FOPC have merits in the following two cases: 

\noindent $\bullet$ The DDSE task has a given computational budget to estimate the state before a new measurement is collected and processed. In this case, in par with the computational demand, FOPC outperforms traditional iterative algorithms.

\noindent $\bullet$ FOPC can perform a prediction of the state while waiting for the measurement to be transmitted from the PMUs to the state estimator; once the measurement is received, the correction step can be performed. The prediction stage is shown to enable substantial improvements in terms of tracking performance.  

Relative traditional approaches based on Kalman filtering, it is worth pointing out that:  (i) FOPC provides an appreciable flexibility to include a variety of performance metrics in the cost function of the problem; for example, the paper will demonstrate how the DDSE can be easily robustified by modifying the cost function. (ii) KF requires covariance matrices; in a DDSE setting where heterogeneous measurements are collected at different rates, it is practically challenging to obtain accurate estimates of noise covariance matrices. And, (iii) KF typically needs an Hessian inverse computations, which can be too computationally burdensome; on the other hand, FOPC relies on first-order updates. Overall, FOPC is naturally data-driven, while KF is grounded on models for the dynamics and the noises. Hence, using the preferred time-varying optimization algorithms, our proposed FOPC-based DDSE methods are more computationally efficient and flexible to incorporate various types of measurements.  


\section{Modeling and Problem Statement}
\label{sec:modeling}


We consider a generic multi-phase unbalanced distribution system with multiphase nodes collected in the set $\cN \cup \{0\}$, $\cN := \{1,\ldots, N\}$, and distribution line segments collected in the set of edges $\cE := \{(m,n)\}$. 
Node $0$ denotes the three-phase slack bus, i.e., the point of connection of the distribution grid with the rest of the electrical system. At each multiphase node, the loads can be either wye- or delta-connected~\cite{Kerstingbook}, with the number of each type being $N^Y$, $N^{\Delta}$ respectively. 

We briefly introduce the AC power-flow model for multiphase distribution systems (a comprehensive description can be found in~\cite{Kerstingbook,bernstein2017load}). To this end, let $\bv$ be a vector collecting the line-to-ground voltages in all phases of the nodes in $\cN$; similarly, vector $\bi$ collects all the phase net current injections, $\bi^{\Delta}$ the phase-to-phase currents in all the delta connections, and vectors $\bs^Y$ and $\bs^\Delta$ collect the net complex powers injected for wye- and delta-connected loads, respectively. Let $N_\phi$ denote the total number of single-phase connections; for example, if all the nodes are three-phase, it follows that $N_\phi = 3N$. 
With these definitions in place, the AC power-flow equations can be compactly written as:   
\begin{subequations} \label{eq:lf}
	\begin{align} 
	& \diag(\bH^\sfT (\bi^{\Delta})^*)\bv + \bs^Y = \diag(\bv) \bi^*, \label{eq:lf_balance}\\
	& \bs^{\Delta} = \diag\left(\bH \bv \right)(\bi^{\Delta})^*, \bi = \bY_{L0} \bv_0 + \bY_{LL} \bv , \label{eq:lf_i} 
	\end{align}
\end{subequations}
where 
$[\bY_{L0}, \bY_{LL}] \in \comps^{N_\phi \times (3+N_\phi)}$ 
is the submatrix of the admittance matrix $\bY$ by eliminating the slack-bus rows,
while $\bH$ is a $N_\phi \times N_\phi$ block-diagonal matrix mapping $\bi^{\Delta}$ to line currents; see~\cite{bernstein2017load} for a detailed description. For real-valued notations, define the $2N_\phi \times 1$ rectangular-form voltage vector
\begin{align} 
\label{eq:voltage}
\bz := [\Re\{\bv\}^\sfT , \Im\{\bv\}^\sfT]^\sfT. 
\end{align}
Similarly, define the vectors $\bu^Y := [\Re\{\bs^Y\}^\sfT , \Im\{\bs^Y\}^\sfT]^\sfT$ and $\bu^\Delta := [\Re\{\bs^\Delta\}^\sfT , \Im\{\bs^\Delta\}^\sfT]^\sfT$ for power variables. 

The proposed approach leverages a linearized AC power flow model to facilitate the development of computationally affordable algorithms that can be implemented in real time. For any given complex powers $\bs^Y, \bs^\Delta$ and its voltage solution $\bv$, this paper leverages a fixed-point approximation of \eqref{eq:lf} as detailed in \cite{bernstein2017load}, to obtain a linearized model that exactly touches upon $(\bs^Y, \bs^\Delta,\bv)$ and $(\mathbf 0, \mathbf0, \bw)$, where $\bw$ is the zero-load voltage solution. Notice that, relative to alternative linearization techniques (e.g.,~\cite{bolognani2015linear,sairaj2015linear,Baran89,sulc2014optimal}), the approach in \cite{bernstein2017load} accounts for both wye and delta connections. 

Accordingly, consider the following linearized model of \eqref{eq:lf}:
	\begin{align} 
    \label{eq:linear_model_v}
	\tilde{\bz}& = \bM^Y \bu^Y + \bM^\Delta \bu^\Delta + \mathbf{m} = \bM \bu  + \mathbf{m}
	\end{align}
where the
model parameters $\bM^Y \in \mathbb{R}^{2N_\phi \times 2N^Y_\phi}$, $\bM^\Delta \in \mathbb{R}^{2N_\phi \times 2N^\Delta_\phi}$, and $\mathbf{m} \in \mathbb{R}^{2N_\phi \times 1}$, with $\bM := [\bM^Y, \bM^\Delta]$ and $\bu := [(\bu^Y)^\sfT, (\bu^\Delta)^\sfT]^\sfT$. The vector $\mathbf{m} := [\Re\{\mathbf{w}\}^\sfT, \Im\{\mathbf{w}\}^\sfT]^\sfT$ ensures that $(\mathbf 0,\mathbf 0, \bw)$ always satisfies the model \eqref{eq:linear_model_v} for $\bu = \mathbf 0$. In addition, matrices $\bM^\Delta$ and $\bM^Y$
are computed based on the network parameters and a given voltage profile. It is worth re-iterating that a linear model is leveraged in this paper to obtain the time-varying convex optimization model described in the next section and, based on that, synthesize dynamic algorithmic solutions with provable tracking properties.    

\subsection{Problem statement}
\label{sec:problem_statement}

We introduce the measurement model for our $\mu$PMU assisted DDSE problem. Assume that the temporal axis is discretized as $t_k = h k$, where $k = 0, 1, \ldots, $ and the sampling period is $h = t_{k+1}-t_{k}$. Assume that $\mu$PMUs are located  at a subset of nodes $\cM_\bv \subset \cN$; let $M_\bv$ denote the number of line-to-ground voltage measurements collected from the multi-phase nodes in $\cM_\bv$. These $\mu$PMUs can obtain accurate measurements of the voltages in rectangular coordinates, and they can produce measurements in real time.  Furthermore, net injected powers from wye and delta connections are measured at the multi-phase nodes in  $\cM_\bu^Y$ and $\cM_\bu^\Delta$, respectively, with $\cM_\bu := \cM_\bu^Y \bigcup \cM_\bu^\Delta \subset \cN$. Accordingly, let $M_\bu$ denote the number of nodes in $\cM_\bu$. 

In this setting, at each time instant $t_k$ a set of new measurements are collected and processed for DSSE; the interval $h$ can be small -- even on the order of seconds -- if the fast-acting measurement capabilities of $\mu$PMUs and distributed energy resources are leveraged. In particular, a (subset of) the following quantities are \emph{measured} at every time $t_k$, $k \in \mathbb{N}$:

\begin{itemize}
	\item $\by_{v}^{(k)} \in \mathbb{R}^{2 M_\bv \times 1}$: measurements of the line-to-ground voltages at all the phases of the nodes $\cM_\bv$. The measurement model for the $\mu$PMUs is $\by_{v}^{(k)}  = \bz_{\cM_\bv}^{(k)} + \bn_{v}^{(k)}$ with  the measurement noise $\bn_{v}^{(k)}$. 
	
\item $\by_{u}^{(k)} \in \mathbb{R}^{2 M_\bu \times 1}$: measurements of the net active and reactive powers from wye and/or delta connections at nodes $\cM_\bu$. The measurement model is $\by_{u}^{(k)} = \bu_{\cM_\bu}^{(k)} + \bn_{u}^{(k)}$ with  the measurement noise $\bn_{u}^{(k)}$. 
	
\end{itemize}

Using~\eqref{eq:linear_model_v}, the measurement equation per time $t_k$ is:   
\begin{align} 
\begin{bmatrix}
\by_{v}^{(k)}  \\
{\by_{u}^{Y}}^{(k)}  \\
{\by^\Delta_{u}}^{(k)} 
\end{bmatrix} 
= &
\begin{bmatrix}
\bJ_{v}\bMY^{(k)} & \bJ_{v}\bMDelta^{(k)} \\
\bJ^Y & \mathbf{0} \\
\mathbf{0} & \bJ^\Delta
\end{bmatrix} 
\begin{bmatrix}
{\bu^Y}^{(k)}  \\
{\bu^\Delta}^{(k)}
\end{bmatrix}
\nonumber \\
& + 
\begin{bmatrix}
\mathbf{m}_v \\
\mathbf{0}  \\
\mathbf{0}  
\end{bmatrix}
+ 
\begin{bmatrix}
\bn_{v}^{(k)}  \\
{\bn^Y_{u}}^{(k)} \\
{\bn^\Delta_{u}}^{(k)}
\end{bmatrix}
\label{eq:meas}
\end{align}
where: $\bJ_{v}$ is a suitable permutation matrix  which selects rows of $\bz$ to form $\bz_{\cM_\bv}$; and, similarly for $\bJ^Y$ and $\bJ^\Delta$,  selecting the measured loads. Recall that $\bMY^{(k)}$  and $\bMDelta^{(k)}$ are time-variant is the last iterate $\bz^{(k)}$ is utilized as the voltage profile for the fixed-point linearization. In our algorithmic development later on, $\bz^{(k)}$ is approximated by its estimated value $\hat{\bz}^{(k)}$ at each time step.

\begin{remark}[Heterogeneous measurements] Although we model the DDSE problem for only voltage phasor data from PMUs and power data from pseudo-measurements or smart meters, it can be generalized to encompass a variety of measurements available in distribution systems, as reviewed in \cite{primadianto2017review}. First, the zero-injection constraint is always satisfied by model \eqref{eq:meas} as only load buses are included by the power variables. Second, as detailed in \cite{bernstein2017load,bernstein2017linear}, the linearized model \eqref{eq:linear_model_v} is generalizable to voltage magnitude, line currents, and power flows. Accordingly, real-time measurements of these variables can be included in \eqref{eq:meas}. To fit various sensing frequency of different measurements, one can set $h$ as the fastest sampling time (typically from PMUs) and maintain the values of other slower measurements until new datum arrives.
\end{remark}
Next, consider rewriting the measurement model~\eqref{eq:meas} in the following compact form: 
\begin{align} 
\label{eq:meas2}
\by^{(k)} = \bG^{(k)} \bu^{(k)} + \bar{\mathbf{m}} + \bn^{(k)}
\end{align}
where $\by^{(k)} := [(\by_{v}^{(k)})^\sfT,({\by_{u}^{Y}}^{(k)})^\sfT, ({\by^\Delta_{u}}^{(k)})^\sfT]^\sfT$, $\bar{\mathbf{m}} := [\mathbf{m}_v^\sfT, \mathbf{0}^\sfT, \mathbf{0}^\sfT]^\sfT$, $\bn^{(k)} := [(\bn_{v}^{(k)})^\sfT,
({\bn^Y_{u}}^{(k)})^\sfT,
({\bn^\Delta_{u}}^{(k)})^\sfT]^\sfT$, and 
\begin{align} 
\bG^{(k)} := \begin{bmatrix}
\bJ_{v}\bMY^{(k)} & \bJ_{v}\bMDelta^{(k)} \\
\bJ^Y & \mathbf{0} \\
\mathbf{0} & \bJ^\Delta
\end{bmatrix} \, .
\end{align}

The data $\mathbf{y}_u$ is usually collected or generated  at much lower quality compared to the high-resolution $\mu$PMU data. In fact, $\mathbf{y}_u$ is either collected from meters or generated from historic load information. Thus, to account for different granularities and precisions, we define the instantaneous  error mismatch loss function $\ell^{(k)}(\bz)$ at time $t_k$ as:
\begin{align} 
\label{eq:costfunction}
\ell^{(k)}(\bu) := &\frac{1}{2}\left\|\by_v^{(k)} - \bG^{(k)}_v\bu - \mathbf{m}_v \right\|_2^2 +  \mathsf{H}\left( {\by^Y_u}^{(k)}-\bJ^Y{\bu^Y} \right) \nonumber \\
& \qquad+ \mathsf{H}\left( {\by^\Delta_u}^{(k)}-\bJ^\Delta{\bu^\Delta}  \right) 
\end{align}
with $\mathbf{G}_v^{(k)} := [\bJ_{v}\bMY^{(k)} , \bJ_{v}\bMDelta^{(k)}]$, and the Huber loss function $\mathsf{H}(\bepsilon) := \sum_{i} \mathsf{H}_i(\epsilon_i)$  written as
\begin{align} \label{eq:huber}
\mathsf{H}_i (\epsilon_i):= \left\{
\begin{array}{ll}
- \delta\epsilon_i-\delta^2/2,& \textrm{if~}\epsilon_i < - \delta \\
|\epsilon_i|^2/2,&\textrm{if~}|\epsilon_i|\leq \delta\\ 
\delta \epsilon_i-\delta^2/2,& \textrm{if~} \epsilon_i > \delta
\end{array}\right.  
\end{align}
where $\delta > 0$ is a positive parameter determined by the load data quality. The Huber loss function is utilized to reject  possible outliers, or down-weight data with substantial measurement errors~\cite{Hastie09}. 

\begin{remark}[Weighted error objective] To accommodate varying data quality,  \textit{non-uniform weights} can be assigned to different types of measurements. With the high accuracy of voltage data from PMUs, one can use a large positive weight for the voltage error term in \eqref{eq:meas}.  Note that this weighted error objective does not affect the problem structure and thus is not included specifically by the algorithmic developments.
\end{remark}

Using the error mismatch $\ell^{(k)}$, 
we formulate the following state estimation problem at time $t_k$: 
\begin{align} 
\label{eq:leastsquares2}
\mathsf{P}^{(k)}(\bu): \min_{\bu \in \mathbb{R}^{N_\phi}} f^{(k)}(\bu):= \ell^{(k)}(\bu) + r^{(k)}(\bu),  k\in \mathbb{N}
\end{align}
where $r^{(k)}(\bu)$ is a (possibly time-varying) regularization function that renders the overall cost function $f^{(k)}(\bu)$ globally strongly convex. 
Notice further that the problem~\eqref{eq:leastsquares2} is unconstrained; however, possible prior information on the minimum and maximum values of the vector $\bu$ can be naturally incorporated in the proposed approach. 

Problem~\eqref{eq:leastsquares2} models a \emph{time-varying} state estimation task under streaming of measurements, and it implicitly defines an \emph{optimal trajectory} $\{\bu^{(k, \star)}\}_{k \in \mathbb{N}}$ for the estimation task. Accordingly, the objective here is to track $\{\bu^{(k, \star)}\}_{k \in \mathbb{N}}$ by processing the incoming measurements in real time. One way to obtain $\{\bu^{(k, \star)}\}_{k \in \mathbb{N}}$ is to solve the problem $\mathsf{P}^{(k)}(\bu)$ \emph{to convergence} (i.e., a batch solution) at each time step $t_k$. However, in a real-time setting with an asynchronous streaming of measurements, a batch solution of~\eqref{eq:leastsquares2} might not be achievable within an interval $h$ due to underlying communication and computational complexity requirements. Thus, the objective of the paper is to develop an algorithmic solution to generate a sequence $\{\bu^{(k)}\}$ of approximate optimizers for the time-varying problem $\{\mathsf{P}^{(k)}(\bu)\}$, which eventually converges to the optimization trajectory $\{\bu^{(k, \star)}\}$. Accordingly, the next section presents a  running  prediction-correction method to solve~\eqref{eq:leastsquares2} in real time. Before doing so, a remark on the strong convexity of $\{\mathsf{P}^{(k)}(\bu)\}$ is in order. 

Strong convexity can ensure that the optimizer trajectory for the sequence of problems $\{\mathsf{P}^{(k)}(\bu)\}$ is unique. In addition, it allows us to establish convergence for the proposed online algorithms.  Regarding the regularization function, possible choices are exemplified next: 

\noindent i) $ r^{(k)}(\bu) = \frac{a}{2}\|\bu - \bu_{\textrm{pr}}^{(k)}\|_2^2$, $a > 0$ and where $\bu_{\textrm{pr}}^{(k)}$ is a priori guess on the load profile; 

\noindent ii) $ r^{(k)}(\bu) = \frac{a}{2}\|\bu\|_2^2$, and it is time invariant;

\noindent iii) if the Hessian of the cost function is available, or, a subspace tracking method is in place, the regularization function can be set to $ r^{(k)}(\bu) = \frac{a}{2}\bu^\sfT \bQ^{(k)} \bu$, where the positive eigenvalues of the matrix $\bQ^{(k)}$ are in the null space of the Hessian of  $\ell^{(k)}(\bu)$. The options i) and ii) would involve a deviation from optimal solutions that one would have obtained by minimizing $\ell^{(k)}(\bu)$, with the magnitude of a possible deviation dependent on the parameter $a$; see~\cite{Koshal11}. The option iii) would not perturb the optimal solution, but the overall solution would incur a higher computational complexity.



\newcommand{\x}{\boldsymbol{x}}

\section{Dynamic state estimation}
\label{sec:method}

We consider \textit{First-Order Prediction-Correction (FOPC)} method~\cite{Simonetto_pc17} to solve the DDSE problem at hand. Inspired by Kalman filtering approaches, prediction-correction approaches allow one to solve a broad class of time-varying convex optimization objectives in a dynamic setting by involving two stages: a prediction phase and a correction phase. In the prediction phase, the algorithm attempts to approach the optimal  solution of the next time period (without new observations) by tapping into intrinsic temporal correlations of the cost function; on the other hand, in the correction phase the predicted vector is corrected using the latest measurement. This mechanism improves the response time to external dynamics and shows a good convergence result when the objective function changes smoothly over time, which is the case in the DDSE setting.

\subsection{FOPC algorithm}
Consider a continuously time-varying unconstrained optimization with  objective $f(\x; t)$ to model external dynamics, along with the sampled counterpart $\{f^{(k)} := f(\x; t_k)\}$ for $t_k = k h$, $k \in \mathbb{N}$. The goal is to produce a trajectory $\{\x^{(k)}\}$ such that $\x^{(k)} \approx \x^*(t_k)$,  where $\x^*(t_k)$ denotes the optimal solution at time $t_k$. In order to predict the solution at time $t_{k+1}$, a strategy is to find  $\x^{(k+1|k)}$ that satisfies the condition
\begin{equation}\label{cond:gamma}
\nabla _x f^{({k+1})}(\x^{(k+1|k)}) = (1-\gamma) \nabla _x f^{({k})}(\x^{(k)}),
\end{equation}
with $\gamma \in [0,1]$.  Varying $\gamma$, this condition imposes optimality ($\gamma = 1$), or the fact that the estimate $\x^{(k+1|k)}$ is no worse than $\x^{(k)}$ in terms of suboptimality ($\gamma = 0$) even when the function changes. The choice $\gamma = 1$ combines moving towards the optimizer while moving with the varying objective function; the choice $\gamma = 0$ represents a \emph{rigid} motion with the objective.  

Notwithstanding the choice of $\gamma$, condition~\eqref{cond:gamma} cannot be computed at time $t_k$ without information about $f^{({k+1})}$. Instead, consider the following Taylor approximation:
\begin{align}
& (1-\gamma) \nabla_x f^{(k)}(\x^{(k)}) =  \nabla_x f^{(k+1)}(\x^{(k+1|k)}) \approx \nabla_x f^{(k)}(\x^{(k)}) \nonumber\\ 
& + \nabla_{xx}f^{(k)}(\x^{(k)})(\x^{(k+1|k)}-\x^{(k)}) + h\nabla_{tx}f^{(k)}(\x^{(k)}) . \label{eq2}
\end{align}
By solving this equation, we have a recursion of the form: 
\begin{multline}\label{gtt-update}
\x^{(k+1|k)} =  \x^{(k)} -\nabla_{xx}f^{(k)}(\x^{(k)})^{-1}\times \\ 
\left(\gamma \nabla_{x}f^{(k)}(\x^{(k)}) + h\nabla_{tx}f^{(k)}(\x^{(k)})\right),
\end{multline}
which (as anticipated) combines a Newton's step with a rigid motion with the objective function.

Concerning over the cost of computing the inverse of the Hessian in a possibly small time interval $h$, the FOPC method further involves a first-order update to solve \eqref{eq2}. Specifically,  the prediction solution of \eqref{eq2} is sought by constructing an equivalent quadratic optimization as follows
\begin{align}
\x^{(k+1|k)} =  \mathrm{argmin}_{\bx} ~~& \hat{f}^{(k)}(\bx)
\end{align}
where 
\begin{align}
\hat{f}^{(k)}(\bx) :=& \frac{1}{2}\bx^\sfT\nabla_{xx}f^{(k)}(\x^{(k)})\bx + \Big(\gamma\nabla_x f^{(k)}(\x^{(k)}) \nonumber \\ 
-& \nabla_{xx}f^{(k)}(\x^{(k)})\x^{(k)} + h\nabla_{tx}f^{(k)}(\x^{(k)})\Big)^\sfT\bx. \label{fopc-quadopt}
\end{align}
Thus, one can replace the update~\eqref{gtt-update} with the gradient descent solution for~\eqref{fopc-quadopt}; i.e., each iteration $p$ is given by: 
\newcommand{\xhat}{\hat{\x}}
\begin{align}
\xhat_{p+1} = & \xhat_p - \alpha \big[\nabla_{xx}f^{(k)}(\x_k)(\xhat_p - \x^{(k)})  + \gamma \nabla_x f^{(k)}(\x_k)\nonumber \\
& + h\nabla_{tx}f^{(k)}(\x_k)\big], ~~~p=0,\ldots,P-1  \label{fopc-gdstep}
\end{align}
where integer $P$ is the number of \textit{prediction steps} that one can afford within an interval $h$, and $\alpha > 0$ the stepsize to be designed later on. Hence, the predicted solution is set to $x^{(k+1|k)} = \xhat_P$. 

Initializing at $x^{(k+1|k)}$, the correction phase further involves $C$ first-order gradient steps; that is, for each $c=0,\ldots,C-1$ 
\begin{align}
\xhat_{c+1} =  \xhat_c - \beta \nabla_{x}f^{(k+1)}(\xhat_{c})(\xhat_c - \x^{(k)}) \label{fopc-gdstepcorrection}
\end{align}
with the stepsize $\beta > 0$. The final corrected estimate is set as $\x^{(k+1)} = \xhat_{C}$. The number of steps $P$ and $C$ for the prediction and correction stages, respectively, are selected based on the duration of the interval $h$.   

The complete running FOPC algorithm for the DDSE problem \eqref{eq:leastsquares2} is tabulated as Algorithm~\ref{alg:pc}. Selection strategies for $P, C$ as well as the stepsizes $\alpha$ and $\beta$ will be elaborated in the next subsection. Notice that, once the prediction $\hat{\bu}^{(k|k-1)}$ and the corrected estimate $\hat{\bu}^{(k)}$ are obtained, the voltage vectors can be readily calculated using the linearized model~\eqref{eq:linear_model_v}. Alternatively, given the estimate of the power injections, the voltages can be calculated by solving the AC power flow equations (see, e.g., the fixed-point power flow method in~\cite{bernstein2017load}).

\begin{algorithm}[h]
\caption{FOPC for DDSE }
\label{alg:pc}

\textbf{Notation}: $\hat{\bu}^{(k)}$ and $\hat{\bu}^{(k|k-1)}$ are the estimate/prediction of $\bu^{(k)}$ at time $t_k$ and $t_{k-1}$ respectively, $k \in \mathbb{N}$.

\textbf{Initialization}:
Choose the number of prediction and correction steps $P,C$, the stepsizes $\alpha, \beta$, the parameter  $\gamma\in[0,1]$, the parameter $\delta$ of the Huber loss, the regularizer $r$. Set $\mathbf{m} = [\Re\{\mathbf{w}\}^\sfT , \Im\{\mathbf{w}\}^\sfT]^\sfT$ and $\bu^{(0)} = \mathbf{0}$.

\textbf{Algorithm}: for $k = 0, 1, 2, \cdots,$ perform: 

\textcolor{black}{At time $t_{k-1}$}:

\emph{Prediction step}: 

[S1-0] When $k=0$, set $\hat{\bu}^{(0|-1)} = \bu_0, \hat{\bz}^{(0|-1)} = \mathbf{m}$ and skip the following prediction.

[S1-1] Set $\bar{\bu}_0 = \hat{\bu}^{(k-1)}$. 

[S1-2] For $p = 0, \ldots, P-1$, do: 
\begin{align} 
\bar{\bu}_{p+1} &:= \bar{\bu}_{p} - \alpha \left[ \left(\nabla_{\bu\bu}f^{(k-1)}(\hat{\bu}^{(k-1)}) \right)(\bar{\bu}_{p} - \hat{\bu}^{(k-1)})   \right. \nonumber \\
&  \left. \hspace{-.2cm} + h \nabla_{t \bu}f^{(k-1)}(\hat{\bu}^{(k-1)}) + \gamma \left(\nabla_{\bu}f^{(k-1)}(\hat{\bu}^{(k-1)}) \right) 
\right]
\label{eq:prediction_step}
\end{align}

[S1-3] Set $\hat{\bu}^{(k|k-1)} = \bar{\bu}_{P}$ and compute \begin{align}
\hat{\bz}^{(k|k-1)} = \begin{bmatrix}
  \bMY^{(k-1)} & \bMDelta^{(k-1)} 
\end{bmatrix}\hat{\bu}^{(k|k-1)} + \mathbf{m}.
\end{align}

\noindent \textcolor{black}{At time $t_{k}$}: 

\noindent \emph{Function update}:

[S2-1] Compute the updated AC linear model $\bMDelta^{(k)}$ and $\bMY^{(k)}$ using the current estimate of voltage, \textit{i.e.}, $\hat{\bz}^{(k|k-1)}$. This is used to update the linear model $\mathbf{G}^{(k)}$ in \eqref{eq:meas2}.

[S2-2] Update $f^{(k)}$ from $f^{(k-1)}$ using $\mathbf{G}^{(k)}$ and the new observations $\by^{(k)}$.

\noindent \emph{Corrections step}: 

[S3-1] Set $\bar{\bu}_0 = \hat{\bu}^{(k|k-1)}$. 

[S3-2] For $c = 0, \ldots, C-1$, do: 
\begin{align} 
\bar{\bu}_{c+1} &:= \bar{\bu}_{c} - \beta  \nabla_{\bu}f^{(k)}(\bar{\bu}_{c}) 
\label{eq:correction_step}
\end{align}

[S3-3] Set $\hat{\bu}^{(k)} = \bar{\bu}_{C}$ and compute 
\begin{align}\hat{\bz}^{(k)} = \begin{bmatrix}
  \bMY^{(k)} & \bMDelta^{(k)} 
\end{bmatrix}\hat{\bu}^{(k)} + \mathbf{m} \, .
\end{align}

\end{algorithm}

\subsection{Online tracking results}

This subsection describes in which sense the sequence of approximate optimizers $\{\hat{\bu}^{(k)}\}_{k\in\mathbb{N}}$ generated by Algorithm~\ref{alg:pc} \emph{tracks} the sampled solution trajectory $\bu^{(*,k)}$. To this end, we will adapt some of the results of~\cite{Simonetto_pc17} to the DDSE problem. 

The following assumptions on the cost function of~\eqref{eq:leastsquares2} and its time variations are presupposed. 

\begin{assumption}\label{as.strong}
The cost function $f^{(k)}(\bu)$ is $\nu$-strongly convex for all $k\in \mathbb{N}$. The regularized function $r^{(k)}(\bu)$ is $L_r$-strongly smooth for all $k\in \mathbb{N}$. 
\end{assumption}

\begin{assumption}\label{as.tv}
The time variation of the gradient of the cost $f^{(k)}(\bu)$
is upper bounded as
\begin{equation}
\|\nabla_{t\bu} f^{(k)}(\bu) \| \leq C_0, \quad \textrm{for all } k \in \mathbb{N}.
\end{equation}
\end{assumption}

Assumption~\ref{as.strong} guarantees that the solution trajectory is unique. This is  why one may need a regularizer $r^{(k)}(\bu)$. As explained in Section \ref{sec:modeling}, this assumption can be verified by utilizing a suitable regularization function. Assumption~\ref{as.tv} makes sure that the functional changes are bounded and thus it is possible to bound the errors arising from a time-varying problem. Note that, in the sequel, the exact value $C_0$ will not be needed for determining the  parameters such as stepsizes, but we will only need the bounded condition. 

First of all, we show that the Hessian of the cost function is lower and upper bounded uniformly in time, as follows. 
 
\begin{proposition} \label{prop:hess}
Under Assumption~\ref{as.strong}, the Hessian $\nabla_{\bu\bu}f^{(k)}(\bu)$ of the cost function~\eqref{eq:leastsquares2} is lower and upper bounded uniformly in time as
\begin{equation}
\nu \leq \|\nabla_{\bu\bu}f^{(k)}(\bu)\| \leq L, \quad \textrm{ for all } k\in \mathbb{N}.
\end{equation}
\end{proposition}

\begin{IEEEproof}
The lower bound follows directly from Assumption~\ref{as.strong}. The upper bound follows by direct computation. The Hessian of $f^{(k)}(\bu)$ is bounded as 
\begin{multline}
\|\nabla_{\bu\bu}f^{(k)}(\bu)\| = \|\nabla_{\bu\bu}\ell^{(k)}(\bu) + \nabla_{\bu\bu}r^{(k)}(\bu)\|\leq \\ \|{\bG^{(k)}_v}^\sfT \bG^{(k)}_v\| + \|{\bJ^Y}^\sfT\bJ^Y\|+ \|{\bJ^\Delta}^\sfT\bJ^\Delta\| + L_r.\nonumber
\end{multline}
By properly defining $L$ as the upper bound of the right-hand term, the Hessian is upper bounded.
\end{IEEEproof}

We are now ready for the online tracking result. 

\begin{theorem}[Adapted from~Theorem 3 of \cite{Simonetto_pc17}]\label{th.1}
Consider the sequence $\{\hat{\bu}^{(k)}\}_{k\in\mathbb{N}}$ generated by Algorithm~\ref{alg:pc}, and let Assumptions~\ref{as.strong}-\ref{as.tv} hold true. Let $\bu^{(*,k)}$ be the optimizer of~\eqref{eq:leastsquares2} at time $t_k$. Choose stepsizes $\alpha$ and $\beta$ as
\begin{equation}\label{cond:step}
\alpha < 2/L, \quad \beta < 2/L,
\end{equation} 
and define the following non-negative quantities 
\begin{equation}
\varrho_{\mathrm{P}} = \max\{|1-\alpha \nu|,|1- \alpha L|\}, \, \varrho_{\mathrm{C}} = \max\{|1-\beta \nu|,|1- \beta L|\}.
\end{equation}
Further, select the number of correction steps $C$ in a way that  
\begin{equation}\label{eq.cvrt_u}
\tau_0 := \varrho_{\mathrm{C}}^{C} \left[\varrho_{\mathrm{P}}^{P} +  (\varrho_{\mathrm{P}}^{P}+1)\Big(1-\gamma +\gamma \frac{2L}{\nu}\Big)\right]  <1.
\end{equation} 
Then, the sequence $\{\|\bu^{(k)}-\bu^{(*,k)}\|\}_{k \in \mathbb{N}}$ converges linearly with rate $\tau_0$ to an asymptotic error bound, and 
\begin{equation}\label{main.result2_0}
\limsup_{k\to\infty}\|\bu^{(k)}-\bu^{(*,k)}\| =O(\varrho_{\mathrm{C}}^{C} h).
\end{equation}
\end{theorem}


Theorem~\ref{th.1} asserts that the sequence $\{\hat{\bu}^{(k)}\}_{k\in\mathbb{N}}$ generated by Algorithm~\ref{alg:pc} \emph{converges to} (and tracks) the sampled solution trajectory $\bu^{(*,k)}$  up to an asymptotic bound. This bound is linearly related to the sampling period $h$ and exponentially decreasing with $C$. Theorem~\ref{th.1} requires standard bounds on  the stepsizes~\eqref{cond:step} and a condition on the number of prediction and correction steps~\eqref{eq.cvrt_u}. 

If one chooses the parameter $\gamma = 0$ (i.e., rigid motion only), then the second condition boils down to
\begin{equation}\label{eq.cvrt_u_g_0}
\varrho_{\mathrm{C}}^{C} \left[2 \varrho_{\mathrm{P}}^{P} +1\right]  <1.
\end{equation}
Since both $\varrho_{\mathrm{P}}$ and $\varrho_{\mathrm{C}}$ are less than $1$ by construction, the condition \eqref{eq.cvrt_u_g_0} is not very restrictive. In fact, fixing a level of prediction $\bar{P}$, then the number of correction steps one has to perform is
\begin{equation}\label{eq.cvrt_u_g_0_inv}
C \geq \left\lceil -\frac{\log(2 \varrho_{\mathrm{P}}^{\bar{P}} +1)}{\log \varrho_{\mathrm{C}}}\right\rceil.
\end{equation}
For reasonable values such as $\varrho_{\mathrm{C}} = \varrho_{\mathrm{P}} = 0.8$ and $\bar{P} = 4$, it implies that $C \geq 3$ is sufficient. When $\gamma > 0$, Condition~\eqref{eq.cvrt_u} gets more restrictive, while the tracking error accuracy may benefit from the Newton's step in \eqref{gtt-update}. 

The choice of $\gamma$ is important to trade-off convergence region, requirements for prediction and correction steps, and conditioning on the measurement matrix. On one hand, if $\gamma = 0$, then the convergence region is bigger, Condition~\eqref{eq.cvrt_u} is less restrictive, which is good when one can afford only a small number of prediction and correction steps (in a fast sampling scenario); however $\gamma = 0$ could be more affected by a badly conditioned measurement matrix. On the other hand, if $\gamma =1$, you get a Newton step in the prediction that helps in case the measurement matrix is badly conditioned, the convergence region gets smaller, and the number of required prediction and correction steps gets higher. 


\section{Numerical Simulations}
\label{sec:simulation}

\begin{figure}[t]
\centering
\includegraphics[width=.9\linewidth]{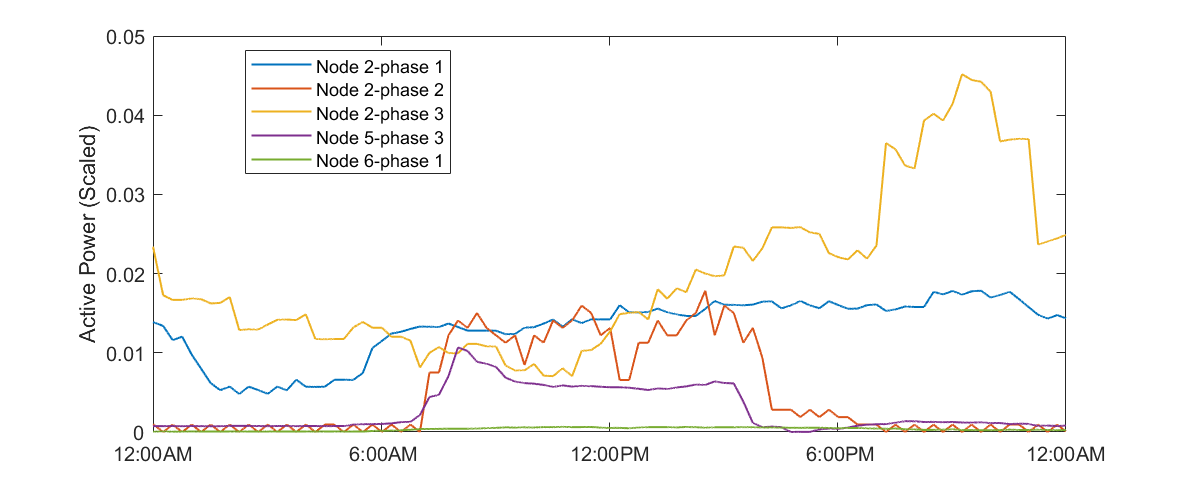}
\caption{Five sample active-power trajectories from the 24-hour load profile dataset.}
\label{fig:power_profile} 
\vspace*{-8mm}
\end{figure}

Numerical tests have been performed using the IEEE 37-bus and 123-bus test feeders (see e.g., \cite{schneider2018analytic} for a description of the feeders) on a standard laptop with Intel$^{\textregistered}$ Core$^{\mathrm{TM}}$ i7-7500 CPU @2.70Hz. The 37-bus test feeder consists of 32 nodes (phases) that are connected to non-zero loads, all delta-connected. The 123-bus test feeder is a popular case with 72 non-zero load nodes and various single-, two-, and three-phase lines, with a mix of delta and wye connections. The load profiles were generated from a real dataset that the National Renewable Energy Laboratory produced from real consumption data received from a utility company in California; the data includes 55 load consumption trajectories over the course of 24-hour, at a time-resolution of 6 seconds. The trajectories for the active power are described in more details in~\cite{FeedbackbasedUnified}; 5 representative trajectories are shown in Fig.~\ref{fig:power_profile}. A constant power factor of 0.95 has been postulated to create the trajectories for the reactive power. The load profiles were randomly chosen for each load node. 

As for the measurement settings, the DDSE algorithm used both the load power data and PMU voltage data. We assume the load active/reactive power injections are observed everywhere in the system. Nonetheless, since they are typically collected by smart meters, the load profiles in Fig. \ref{fig:power_profile} were downsampled to generate the measurement data at a slower time resolution of 10 minutes. Specifically, the power measurements were kept as the average value within every 10-minute window. Accordingly, the parameter for the Huber's loss function in \eqref{eq:huber} was chosen to be $\delta = 8{e}{-4}$, coinciding with the deviation level between the down-sampled power measurement and the actual value.  Since PMUs are not typically installed everywhere, only selected buses are assumed to be equipped with high-quality voltage measurements at the same resolution as the load data. Addictive Gaussian noises with a small standard deviation of $1{e}{-5}$ were added to the actual voltage profiles to reflect the high sensing capability of PMUs. Accordingly, we set the weights of voltage measurements to be $1{e}{3}$ to match the small noise variation therein. 

\subsection{37-Bus Test Feeder}

This test feeder was used to compare the FOPC-based DDSE algorithm under different computational settings, with the performance in terms of both the tracking error with the instantaneous optimum and estimation error with the actual voltage state. For simplicity of implementation, the regularization term in \eqref{eq:leastsquares2} was chosen to be $r(\mathbf u) = \frac{1}{2} \|\mathbf u\|_2$. This sufficiently small term guarantees convergence and yet does not degrade significantly the performance of the estimation problem (as we will show).  To enhance the effectiveness of the prediction step under potentially ill-conditioning issue (see Sec. \ref{sec:method}), we set the parameter $\gamma =0.9$ to be close to 1.

Based on these parameter settings, the norm of the Hessian matrix as in Prop. \ref{prop:hess} can be bounded within the interval $[3.5e4,~4e4]$. Thus, the stepsize parameters $\alpha$ and $\beta$ were chosen  as $1e{-4}$ according to Theorem \ref{th.1}. Accordingly, the two parameters therein, $\varrho_{\mathrm{P}}$ and $\varrho_{\mathrm{C}}$, are around $0.65$, 
both less than 1. To satisfy the convergence condition in \eqref{eq.cvrt_u}, it is sufficient to have $C=5$ steps of correction even if no prediction is performed ($P=0$), as $\tau_0 = 0.8$ in this case.

\subsubsection{Fixed number of correction steps}
We first show the advantage to add prediction phase before new measurements are processed. As mentioned in Sec. \ref{sec:intro}, FOPC can perform a prediction of the state while waiting for the measurement to be transmitted from the PMUs; once the measurement is received, a fixed number of correction steps can be performed. Thus, this test uses a fixed $C=5$ steps of correction and compares the results from $P=0,5,10$ steps of prediction. 

\begin{figure}[t]
\centering
\includegraphics[width=.8\linewidth]{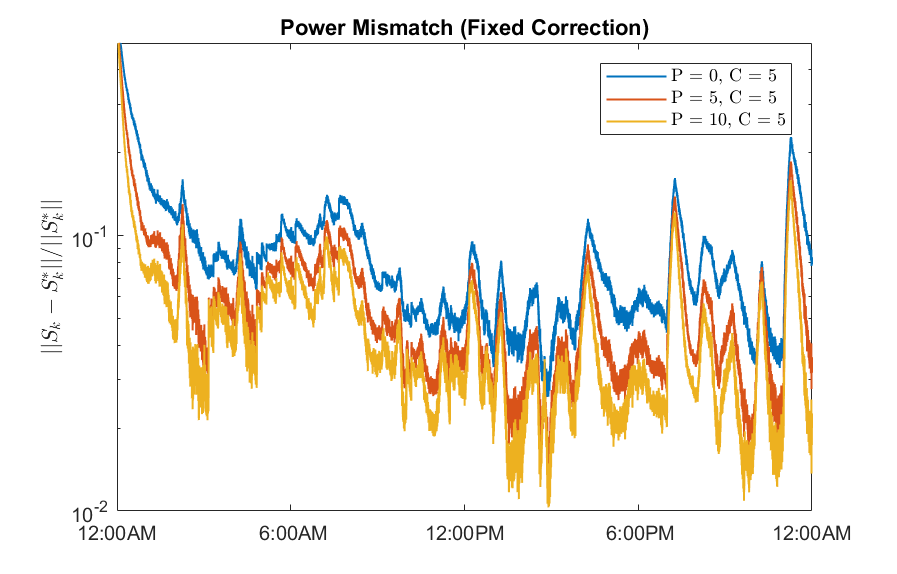}
\caption{Relative error of tracking the instantaneous optimal system-wide power state under fixed $C$ steps. 
}\label{fig:Fixed-Correction}
\vspace*{-6mm}
\end{figure}

\begin{figure}[t]
\centering
\includegraphics[width=.8\linewidth]{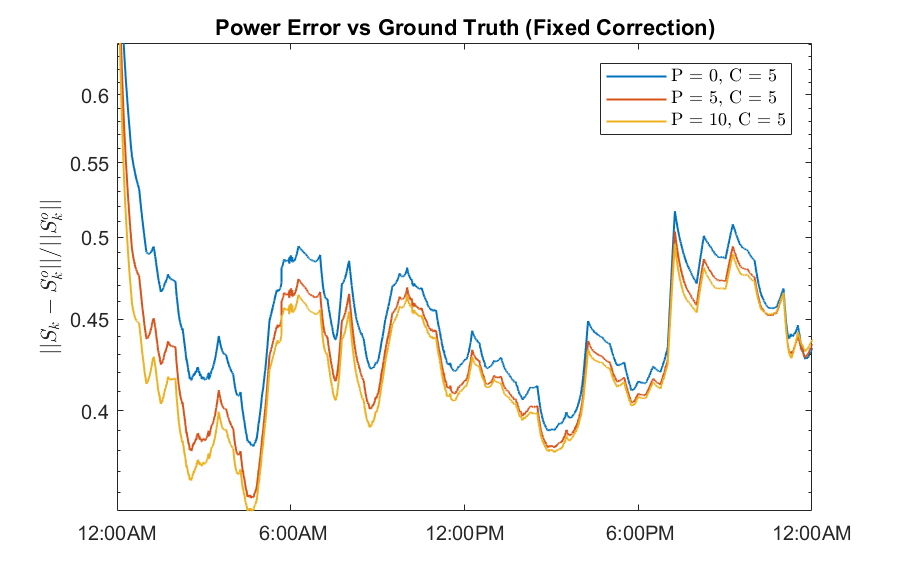}
\includegraphics[width=.8\linewidth]{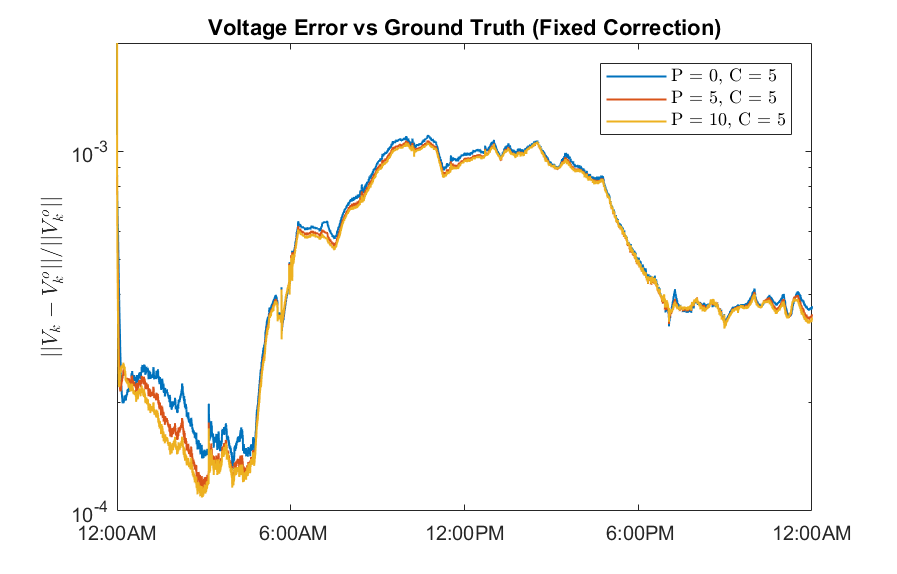}
\caption{Relative error of estimating the ground-truth system-wide power state (top) and voltage output (bottom) under fixed $C$ steps.} \label{fig:Fixed-Correction-Truth} 
\vspace*{-5mm}
\end{figure}

Fig.~\ref{fig:Fixed-Correction} plots the relative tracking error with the instantaneous optimal state $\bu^{(*,k)}$. 
Clearly, a larger number of $P$ does help the tracking of the optimal solution, with more noticeable change in mismatch error from $P=0$ to $P=5$ steps of prediction. It also makes the relative error trajectory more quickly to reach the steady-state level of below $0.1$. Hence, the prediction phase has been shown to improve the tracking error performance for the time-varying DDSE problem.

Furthermore, we compare the relative estimation errors for both the system-wide power states and corresponding voltage outputs as compared to the ground-truth values, as plotted in Fig. \ref{fig:Fixed-Correction-Truth}. As for the power states, the relative estimation error increases from the tracking error level. This could be due to the approximation error of the linearized model adopted by the objective function and the measurement error from low-resolution power data. Meanwhile, the effectiveness of prediction phase is still evident in improving the estimation error and convergence rate. More interestingly, the voltage estimation error is very minimal at the level of below 1e-3 for all scenarios, corroborating that the PMU voltage data is instrumental for recovering the feeder voltage profile even under highly uncertain power measurements. 

\subsubsection{Fixed computational time}
This test compares the FOPC performance under a total computational time constraint.  
The computational time in the prediction phase is mainly spent on a one-shot computation of Hessian $\nabla_{\bu\bu}f^{(k)}(\hat{\bu}^{(k)})$ in~\eqref{eq:prediction_step}, while that of the correction phase grows linearly with $C$ as it needs to compute the gradient in every step. 
We choose two sets of $(P,C)$ values: $(8,3)$ and $(0,6)$, both taking roughly a total of $0.3$ms per iteration.
For the case of $(P,C)=(8,3)$, we verified again the convergence condition in \eqref{eq.cvrt_u} with $\tau_0 \approx 0.8$. 

\begin{figure}[t]
\centering
\includegraphics[width=.8\linewidth]{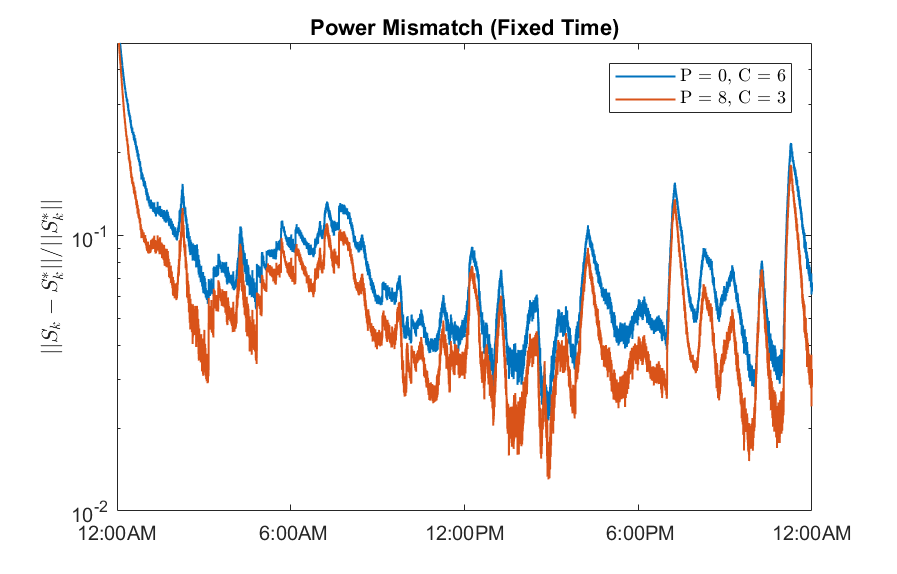}
\caption{Relative error of tracking the instantaneous optimal system-wide power state under a fixed computational time.}
\label{fig:Fixed-Time}
\vspace*{-4mm}
\end{figure}

\begin{figure}[t]
\centering
\includegraphics[width=.8\linewidth]{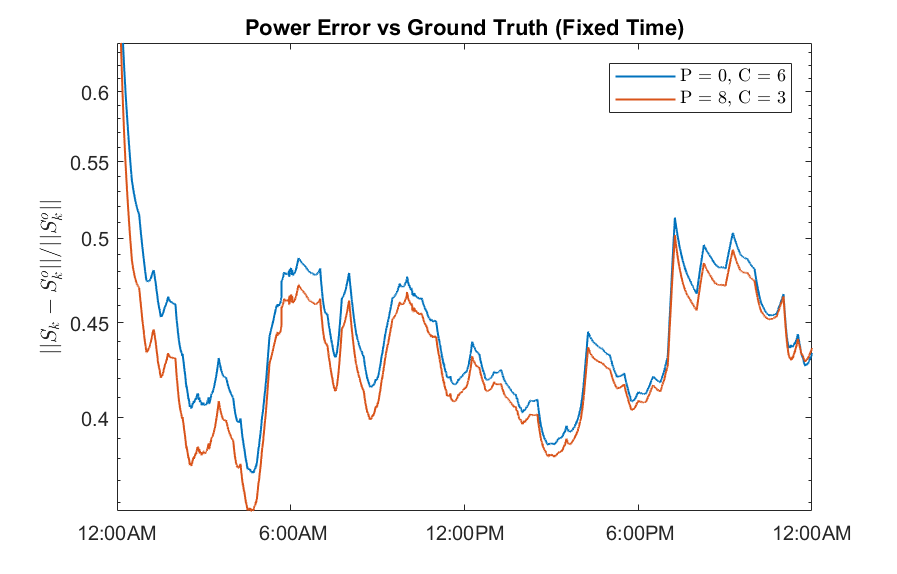}
\includegraphics[width=.8\linewidth]{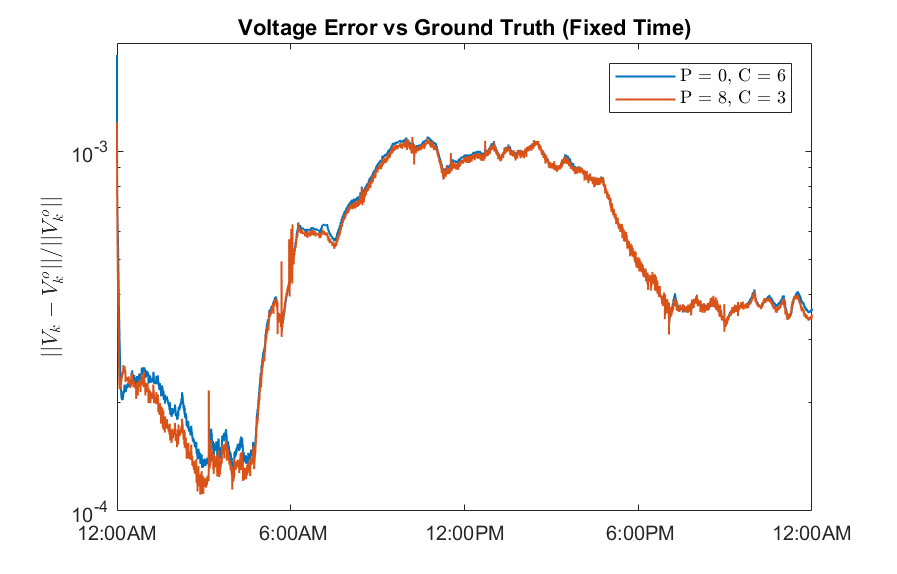}
\caption{Relative error of estimating the ground-truth system-wide power state (top) and voltage output (bottom) under a fixed computational time.}
\label{fig:Fixed-Time-Truth}
\vspace*{-5mm}
\end{figure}

Fig.~\ref{fig:Fixed-Time} and Fig.~\ref{fig:Fixed-Time-Truth} plot the relative tracking and estimation error trajectories for the scenarios, respectively, as in the last test. Interestingly, even with smaller $C$ to compensate for the one-shot Hessian computation, the improvement of including  the prediction phase can be demonstrated. Both the tracking and estimation errors for the power states are lower for the case of $(P,C)=(8,3)$. The voltage error is again very small thanks to the high-quality PMU data. 

\subsubsection{Varying number of PMUs} 
\begin{figure}[t]
\centering
\includegraphics[width=.72\linewidth]{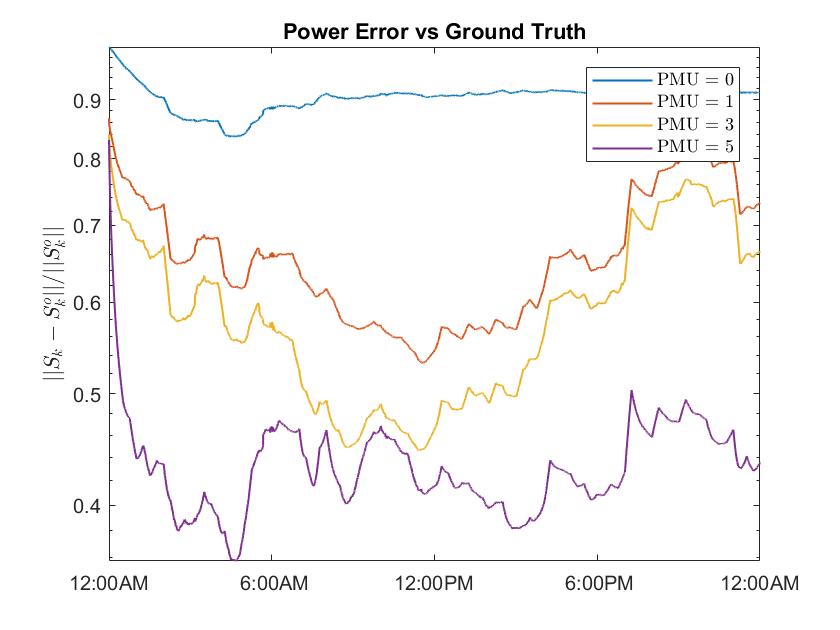}
\includegraphics[width=.72\linewidth]{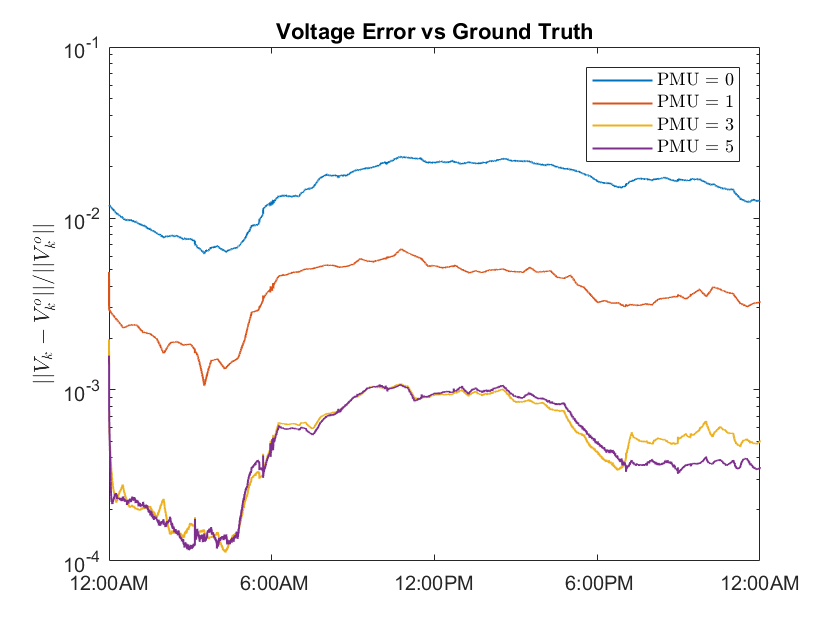}
\caption{Relative error of estimating the ground-truth system-wide power state (top) and voltage output (bottom) under varying number of PMUs.}
\label{fig:PMU-Comparison-Truth}
\vspace*{-5mm}
\end{figure}

We further compare the performance when different number of PMUs are installed in the system.
Fig.~\ref{fig:PMU-Comparison-Truth} plots the relative error of estimating both the ground-truth power and voltage variables. Clearly, more high-quality voltage data can significantly improve the estimation error performance. For the 37-bus case, it seems that 3 PMUs are sufficient for estimating the system-wide voltage, as there is no noticeable improvement with 5 PMUs. This is not necessarily the case when estimating the power state, due to the high uncertainty of power measurements. 

\subsection{123-Bus Test Feeder}
\begin{figure}[t]
\centering
\includegraphics[width=0.8\linewidth]{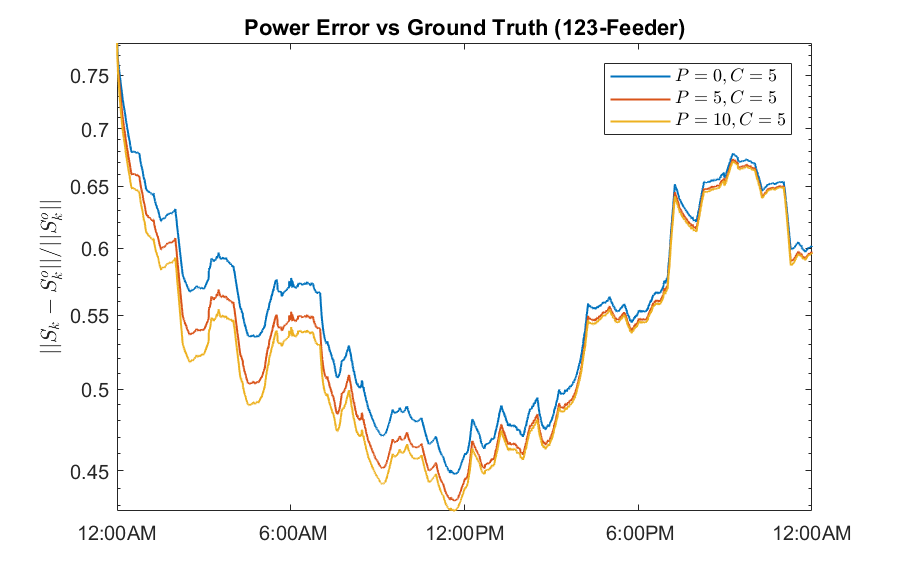}
\includegraphics[width=0.8\linewidth]{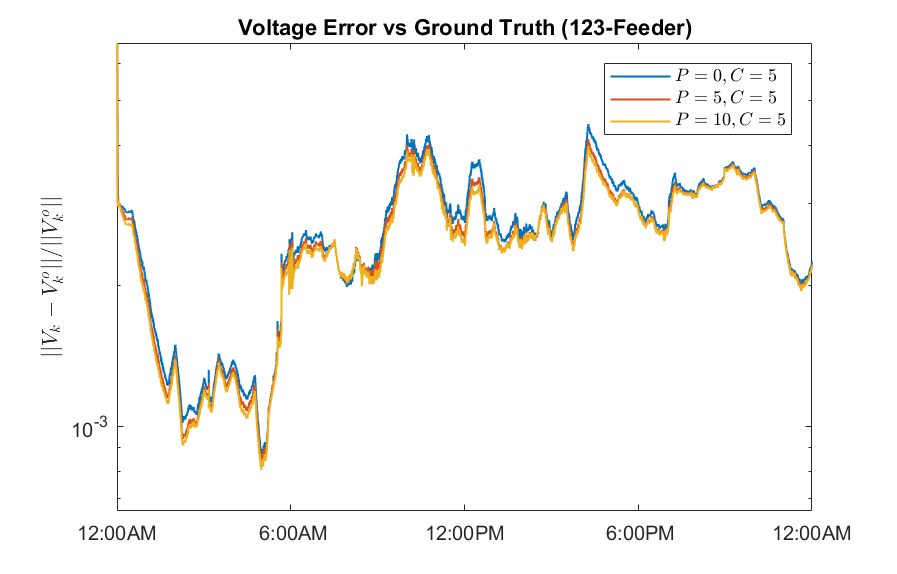}
\caption{Relative error of estimating the ground-truth system-wide power state (top) and voltage output (bottom) of the 123-bus case for a fixed $C$.}
\label{fig:Fixed-Correction-Truth-123}
\vspace*{-4mm}
\end{figure}

Last, we tested the FOPC method on the 123-bus case to demonstrate its scalability. The parameter settings follow from those in the 37-bus tests, with the estimation error comparisons for fixed $C$ given in Fig.~\ref{fig:Fixed-Correction-Truth-123}.  Similar results have been observed, corroborating the improvement of including the prediction phase and high-quality voltage data from PMUs.  

\section{Conclusions}
\label{sec:conc}

This paper presented a distribution state estimation algorithm that can dynamically incorporate fast and accurate PMU voltage data. The first-order prediction-correction (FOPC) algorithm is proposed to solve the time-varying optimization problem of DDSE using a linearized power flow model. Compared to existing recursive updates, the FOPC iterations are computationally simple and require no specific modeling of system transition, suitable for the time-critical DDSE problem where the load dynamics is difficult to model. Numerical tests have shown that the data-driven prediction phase of FOPC is effective in reducing the mismatch error in tracking the power state variable. With the availability of high-quality voltage data, the voltage estimation error is significantly small. Future work includes exploring more diverse types of distribution system measurements and large-scale system validations using real data.

\bibliographystyle{IEEEtran}
\bibliography{biblio,biblio2,DSSEbib}

\begin{thebibliography}{10}
\providecommand{\url}[1]{#1}
\csname url@samestyle\endcsname
\providecommand{\newblock}{\relax}
\providecommand{\bibinfo}[2]{#2}
\providecommand{\BIBentrySTDinterwordspacing}{\spaceskip=0pt\relax}
\providecommand{\BIBentryALTinterwordstretchfactor}{4}
\providecommand{\BIBentryALTinterwordspacing}{\spaceskip=\fontdimen2\font plus
\BIBentryALTinterwordstretchfactor\fontdimen3\font minus
  \fontdimen4\font\relax}
\providecommand{\BIBforeignlanguage}[2]{{%
\expandafter\ifx\csname l@#1\endcsname\relax
\typeout{** WARNING: IEEEtran.bst: No hyphenation pattern has been}%
\typeout{** loaded for the language `#1'. Using the pattern for}%
\typeout{** the default language instead.}%
\else
\language=\csname l@#1\endcsname
\fi
#2}}
\providecommand{\BIBdecl}{\relax}
\BIBdecl

\bibitem{Meier}
A.~von Meier, D.~Culler, A.~McEachern, and R.~Arghandeh, ``Micro-synchrophasors
  for distribution systems,'' in \emph{Proc. IEEE PES Innovative Smart Grid
  Tech. Conf.}, Feb 2014, pp. 1--5.

\bibitem{PMUGomez}
A.~Gomez-Exposito, A.~Abur, P.~Rousseaux, A.~de~la Villa~Jaen, and
  C.~Gomez-Quiles, ``On the use of {PMU}s in power system state estimation,''
  in \emph{17th Power Systems Computation Conf.}, vol.~22, 2011.

\bibitem{primadianto2017review}
A.~Primadianto and C.-N. Lu, ``A review on distribution system state
  estimation,'' \emph{IEEE Transactions on Power Systems}, vol.~32, no.~5, pp.
  3875--3883, 2017.

\bibitem{baran2}
M.~E. Baran and A.~W. Kelley, ``State estimation for real-time monitoring of
  distribution systems,'' \emph{IEEE Trans. Power Syst.}, vol.~9, no.~3, pp.
  1601--1609, Aug 1994.

\bibitem{lu}
C.~N. Lu, J.~H. Teng, and W.~H.~E. Liu, ``Distribution system state
  estimation,'' \emph{IEEE Trans. Power Syst.}, vol.~10, no.~1, pp. 229--240,
  Feb 1995.

\bibitem{baran}
M.~E. Baran and A.~W. Kelley, ``A branch-current-based state estimation method
  for distribution systems,'' \emph{IEEE Trans. Power Syst.}, vol.~10, no.~1,
  pp. 483--491, Feb 1995.

\bibitem{lin2}
W.~M. Lin and J.~H. Teng, ``State estimation for distribution systems with
  zero-injection constraints,'' \emph{IEEE Trans. Power Syst.}, vol.~11, no.~1,
  pp. 518--524, Feb 1996.

\bibitem{EqConKorres}
G.~N. Korres, ``A robust algorithm for power system state estimation with
  equality constraints,'' \emph{IEEE Trans. Power Systems}, vol.~25, no.~3, pp.
  1531--1541, Aug 2010.

\bibitem{HaughtonHeydt}
D.~A. Haughton and G.~T. Heydt, ``A linear state estimation formulation for
  smart distribution systems,'' \emph{IEEE Trans. Power Syst.}, vol.~28, no.~2,
  pp. 1187--1195, May 2013.

\bibitem{LPFMeier}
H.~Ahmadi, J.~R. Mart\'i, and A.~von Meier, ``A linear power flow formulation
  for three-phase distribution systems,'' \emph{IEEE Trans. Power Syst.},
  vol.~31, no.~6, pp. 5012 -- 5021, 2016.

\bibitem{DSSEmuscas}
C.~Muscas, M.~Pau, P.~A. Pegoraro, and S.~Sulis, ``Uncertainty of voltage
  profile in {PMU}-based distribution system state estimation,'' \emph{IEEE
  Trans. Instrum. Meas.}, vol.~65, no.~5, pp. 988--998, May 2016.

\bibitem{sarri2012state}
S.~Sarri, M.~Paolone, R.~Cherkaoui, A.~Borghetti, F.~Napolitano, and C.~A.
  Nucci, ``State estimation of active distribution networks: comparison between
  wls and iterated kalman-filter algorithm integrating pmus,'' in \emph{Proc.
  IEEE PES ISGT Europe}, 2012.

\bibitem{carquex2018state}
C.~Carquex, C.~Rosenberg, and K.~Bhattacharya, ``State estimation in power
  distribution systems based on ensemble kalman filtering,'' \emph{IEEE
  Transactions on Power Systems}, 2018.

\bibitem{zhao2017robust}
J.~Zhao and L.~Mili, ``Robust unscented kalman filter for power system dynamic
  state estimation with unknown noise statistics,'' \emph{IEEE Transactions on
  Smart Grid}, 2017.

\bibitem{valverde2011unscented}
G.~Valverde and V.~Terzija, ``Unscented kalman filter for power system dynamic
  state estimation,'' \emph{IET Generation, Transmission \& Distribution},
  vol.~5, no.~1, pp. 29--37, 2011.

\bibitem{bernstein2017load}
A.~Bernstein, C.~Wang, E.~Dall'Anese, J.-Y. Le~Boudec, and C.~Zhao, ``Load-flow
  in multiphase distribution networks: Existence, uniqueness, non-singularity,
  and linear models,'' \emph{IEEE Trans. Power Sys.}, 2018.

\bibitem{Simonetto_pc15}
A.~Simonetto, A.~Mokhtari, A.~Koppel, G.~Leus, and A.~Ribeiro, ``A class of
  prediction-correction methods for time-varying convex optimization,''
  \emph{IEEE Trans. Signal Proc.}, vol.~64, no.~17, pp. 4576--4591, Sept. 2016.

\bibitem{Simonetto_pc17}
A.~Simonetto and E.~Dall'Anese, ``Prediction-correction algorithms for
  time-varying constrained optimization,'' \emph{IEEE Trans. on Signal Proc.},
  vol.~65, no.~20, pp. 5481--5494, Oct. 2017.

\bibitem{Kerstingbook}
W.~H. Kersting, \emph{Distribution System Modeling and Analysis}.\hskip 1em
  plus 0.5em minus 0.4em\relax 2nd ed., Boca Raton, {FL}: {CRC} Press, 2007.

\bibitem{bolognani2015linear}
S.~Bolognani and F.~D\"orfler, ``Fast power system analysis via implicit
  linearization of the power flow manifold,'' \emph{Allerton Conf.
  Communication, Control, and Computing}, 2015.

\bibitem{sairaj2015linear}
S.~Dhople, S.~Guggilam, and Y.~Chen, ``Linear approximations to {AC} power flow
  in rectangular coordinates,'' \emph{Allerton Conf. Communication, Control,
  and Computing}, 2015.

\bibitem{Baran89}
M.~E. Baran and F.~F. Wu, ``Network reconfiguration in distribution systems for
  loss reduction and load balancing,'' \emph{{IEEE} Trans. on Power Delivery},
  vol.~4, no.~2, pp. 1401--1407, Apr. 1989.

\bibitem{sulc2014optimal}
P.~Sulc, S.~Backhaus, and M.~Chertkov, ``Optimal distributed control of
  reactive power via the alternating direction method of multipliers,''
  \emph{IEEE. Trans. Energy Conversion}, vol.~29, no.~4, pp. 968--977, 2014.

\bibitem{bernstein2017linear}
A.~Bernstein and E.~Dall'Anese, ``Linear power-flow models in multiphase
  distribution networks,'' in \emph{IEEE PES ISGT-Europe}, 2017.

\bibitem{Hastie09}
T.~Hastie, R.~Tibshirani, and J.~Friedman, \emph{{The Elements of Statistical
  Learning: Data Mining, Inference, and Prediction}}.\hskip 1em plus 0.5em
  minus 0.4em\relax Springer, 2009.

\bibitem{Koshal11}
J.~Koshal, A.~Nedi\'{c}, and U.~Y. Shanbhag, ``Multiuser optimization:
  Distributed algorithms and error analysis,'' \emph{{SIAM} J. on
  Optimization}, vol.~21, no.~3, pp. 1046--1081, 2011.

\bibitem{schneider2018analytic}
K.~Schneider, B.~Mather, B.~C. Pal, C.-W. Ten, G.~Shirek, H.~Zhu, J.~Fuller
  \emph{et~al.}, ``Analytic considerations and design basis for the ieee
  distribution test feeders,'' \emph{IEEE Trans. Power Systems}, vol.~33,
  no.~3, pp. 3181--3188, 2018.

\bibitem{FeedbackbasedUnified}
A.~Bernstein and E.~Dall'Anese, ``Real-time feedback-based optimization of
  distribution grids: A unified approach,'' 2017, [Online] Available at:
  https://arxiv.org/abs/1711.01627.

\end{thebibliography}
\end{document}